\documentclass{article}
\usepackage{amsfonts}
\usepackage{amsmath}

\input{tcilatex}
\begin{document}

\begin{center}
\bigskip {\LARGE On Inverse Surfaces in Euclidean 3-Space}

\bigskip

\textbf{\ M. Evren AYDIN}$\dagger $\textbf{, Mahmut ERG\"{U}T}

\bigskip

\textbf{Department of Mathematics, Firat University}

\bigskip

\textbf{Elazig, 23119, Turkey}

\bigskip

$\dagger $\textbf{email:} \textbf{aydnevren@gmail.com}\\[0pt]
\end{center}

\textbf{Abstract. }In this paper, we study the inverse surfaces\ in
3-dimensional Euclidean space $\mathbb{E}^{3}$. We obtain some results
relating Christoffel symbols, the normal curvatures, the shape operators and
the third fundamental forms of the inverse surfaces.\newline

\linebreak Christoffel symbols.\newline

\textbf{Msc.} 11A25, 53A04, 53A05.\newline

\textbf{1. Introduction}\newline

The inversive geometry have been studied for many years by
various\linebreak\ mathematicians. See for instance, $\cite%
{2,4,5,6,9,10,11,13}.$ In $\cite{3,8,12},$ the authors gave inversions of
curves and surfaces. In $\cite{1},$ he obtained some results relating
inverse surfaces.

During the same time that projective geometry was emerging, \linebreak
mathematicians mostly began to deal with circles. Their key tool was
\linebreak inversion. This interest gave rise to a new geometry, called
inversive \linebreak geometry, which was able to provide particularly
striking proofs of previously known results in Euclidean geometry as well as
new results.

This geometry have a wide application area in physics. For example, in the
study of electrostatics let $l_{1}$ and $l_{2}$ be two infinitely long
parallel cylinders of\linebreak\ opposite charge. Then intersection of the
surface of equipotential with a\linebreak\ horizontal plane is two families
of circles (and a single line), and a point charge placed in the
electrostatic field moves along a circular path through a \linebreak
specific point inside each cylinder, at right angles to circles in the
families. This \linebreak geometry deal with such families of circles.

Now let put an inversion. Let $S$ be a sphere with center $c$ and radius $r$%
. Suppose that point $p$ is other than $c$, then the inverse of $p$ with
respect to the sphere $S$ is the point $q$ on line $cp$ such that 
\begin{equation}
\left\Vert \overrightarrow{cp}\right\Vert \left\Vert \overrightarrow{cq}%
\right\Vert =r^{2}.  \tag{1.1}
\end{equation}%
The sphere $S$ is called the sphere of inversion, the point $c$ is the
center of inversion and the radius $r$ is the radius of inversion.

The inversion maps the inside of sphere into outside vice versa and it
\linebreak unchange points on the sphere. $c$ has no image, and no point of
the space is mapped to $c$. But, the points close to $c$ are mapped to
points far from $c$ and vice versa. So, let us say that the image of $c$
"point at infinity" and vice versa.

In this study is organised as follows: Initially, it was obtained that
regularity and singularity for inverse surfaces are properties remaining
unchanged under an inversion. Then, some relations\ concerning the the
normal curvatures, the shape operators, Christoffel symbols, the asymtotic
and the conjugate vectors of inverse surfaces were found.\newline

\textbf{2. Basic Notions}\newline

Let $c$ $\in $ $\mathbb{E}^{n}$ and $r$ $\in 
\mathbb{R}
^{+}$. We denote that $\mathbb{E}^{n^{\ast }}=$ $\mathbb{E}^{n}-\left\{
c\right\} .$ Then, an inversion of $\mathbb{E}^{n}$ with the center $c\in 
\mathbb{E}^{n}$ and the radius $r$ \ is the map 
\begin{equation*}
inversion\left[ c,r\right] :\mathbb{E}^{n^{\ast }}\longrightarrow \mathbb{E}%
^{n^{\ast }}
\end{equation*}%
given by 
\begin{equation}
inversion\left[ c,r\right] \left( p\right) =c+\frac{r^{2}}{\left\Vert
p-c\right\Vert ^{2}}\left( p-c\right) .  \tag{2.1}
\end{equation}%
We use $inv\left[ c,r\right] $ instead of $inversion\left[ c,r\right] $ for
shortness.\newline

The inversion of a curve $\alpha :\left( a,b\right) \longrightarrow $ $%
\mathbb{E}^{n^{\ast }}$ is the curve%
\begin{equation*}
t\longrightarrow inversioncurve\left[ c,r,\alpha \right] \left( t\right) 
\end{equation*}%
defined by 
\begin{equation}
inversioncurve\left[ c,r,\alpha \right] \left( t\right) =inv\left[ c,r\right]
\left( \alpha \left( t\right) \right) .  \tag{2.2}
\end{equation}%
\qquad\ 
\begin{gather*}
\FRAME{itbpF}{1.6604in}{1.6604in}{0in}{}{}{Figure}{\special{language
"Scientific Word";type "GRAPHIC";maintain-aspect-ratio TRUE;display
"USEDEF";valid_file "T";width 1.6604in;height 1.6604in;depth
0in;original-width 2.0176in;original-height 2.0176in;cropleft "0";croptop
"1";cropright "1";cropbottom "0";tempfilename '../Kitap ve
Makaleler/MAKALE/makale 1/BÝRÝNCÝ MAKALE/LCE7FP00.wmf';tempfile-properties
"XPR";}}\text{ \ \FRAME{itbpF}{1.657in}{1.6492in}{0in}{}{}{Figure}{\special%
{language "Scientific Word";type "GRAPHIC";maintain-aspect-ratio
TRUE;display "USEDEF";valid_file "T";width 1.657in;height 1.6492in;depth
0in;original-width 1.8801in;original-height 1.8697in;cropleft "0";croptop
"1";cropright "1";cropbottom "0";tempfilename '../Kitap ve
Makaleler/MAKALE/makale 1/BÝRÝNCÝ MAKALE/LCE7UV02.wmf';tempfile-properties
"XPR";}}} \\
\text{Figure 1: Astroid and its inversion curve with respect to unit circle.}
\end{gather*}%
As shown in Figure 1, the points of astroid which are exterior to circle of
\linebreak inversion are mapped to points on the interior circle of
inversion, and vice versa.\newpage 

Let $x:U\longrightarrow \mathbb{E}^{n^{\ast }}$ be a patch. The inverse
patch of $x$ with respect to $inv\left[ c,r\right] $ is the patch given by

\begin{equation}
y=inv\left[ c,r\right] \circ x.  \tag{2.3}
\end{equation}

\begin{gather*}
\FRAME{itbpF}{1.3197in}{2.6982in}{-0.0208in}{}{}{Figure}{\special{language
"Scientific Word";type "GRAPHIC";display "USEDEF";valid_file "T";width
1.3197in;height 2.6982in;depth -0.0208in;original-width
1.3776in;original-height 4.5377in;cropleft "0";croptop "1";cropright
"1";cropbottom "0";tempfilename '../Kitap ve Makaleler/MAKALE/ÝKÝNCÝ
MAKALE/L5UTNM0I.wmf';tempfile-properties "XPR";}}\text{ \ \ \ \ \ \ \ \ }%
\FRAME{itbpF}{1.2228in}{2.5806in}{0in}{}{}{Figure}{\special{language
"Scientific Word";type "GRAPHIC";display "USEDEF";valid_file "T";width
1.2228in;height 2.5806in;depth 0in;original-width 2.0176in;original-height
3.9565in;cropleft "0";croptop "1";cropright "1";cropbottom "0";tempfilename
'../Kitap ve Makaleler/MAKALE/ÝKÝNCÝ
MAKALE/L5UTNM0J.wmf';tempfile-properties "XPR";}} \\
\text{Figure 2: Circular cylinder and its inverse patch.}
\end{gather*}%
As shown above the figure, the points of the cylinder which are exterior to
the sphere of inversion are mapped to the points on interior, and vice
versa. Note that as the cylinder arise in direction $\pm $ $z$ its inverse
close to center of inversion (origin).

Throughout this paper, we assume that $\Phi $ is $inv\left[ c,r\right] $ and 
\begin{equation*}
y=\Phi \circ x.
\end{equation*}%
Here $x:U\longrightarrow \mathbb{E}^{3^{\ast }}$ is a patch in $\mathbb{E}%
^{3^{\ast }}$.\newline

\textbf{Definition 2.1.} $\left( \cite{7}\right) $ Let $\Phi =inv\left[ c,r%
\right] $. Then, the tangent map of $\Phi $ at $p$ is the map 
\begin{equation*}
\Phi _{\ast p}=T_{p}\left( \mathbb{E}^{n^{\ast }}\right) \longrightarrow
T_{\Phi \left( p\right) }\left( \mathbb{E}^{n^{\ast }}\right)
\end{equation*}%
\ given by%
\begin{equation}
\Phi _{\ast p}\left( v_{p}\right) =\frac{r^{2}v_{p}}{\left\Vert
p-c\right\Vert ^{2}}-\frac{2r^{2}\left\langle \left( p-c\right)
,v_{p}\right\rangle }{\left\Vert p-c\right\Vert ^{4}}\left( p-c\right) , 
\tag{2.4}
\end{equation}%
where $v_{p}\in T$ $_{p}\left( \mathbb{E}^{n^{\ast }}\right) .$\newline

\textbf{Lemma 2.2. }$\left( \cite{7}\right) $ An inversion is a conformal
diffeomorphism.\newline

\textbf{3. The Shape Operators and Christoffel Symbols of Inverse }%
\linebreak \textbf{Surfaces}\newline

Let $x:U\longrightarrow \mathbb{E}^{3^{\ast }}$ be a patch and $\left(
u_{0},v_{0}\right) \in U.$ $x$ is regular at $\left( u_{0},v_{0}\right) $ if 
$x_{u}\times x_{v}$ is nonzero at $\left( u_{0},v_{0}\right) .$ Thus we have
following theorem.\newline

\textbf{Theorem 3.1. \ }Regularity and singularity\ for the inverse surfaces
are invariant properties under the inversion.\newline

\textbf{Proof. }We assume that $\Phi =inv[c,r].$ Then, from $\cite{8}$ we
have

\begin{equation}
y_{u}\times y_{v}=\lambda ^{2}\left( -x_{u}\times x_{v}+\frac{2\left\langle
\left( x-c\right) ,\left( x_{u}\times x_{v}\right) \right\rangle }{%
\left\Vert x-c\right\Vert ^{2}}\left( x-c\right) \right) .  \tag{3.1}
\end{equation}

If $p=x\left( u_{0},v_{0}\right) $ a regular point of $x$, then 
\begin{equation*}
x_{u}\times x_{v}\mid _{p}\neq 0.
\end{equation*}%
Assume that%
\begin{equation*}
y_{u}\times y_{v}\mid _{q}=0,
\end{equation*}%
where $\Phi \left( p\right) =q\in y\left( U\right) .$ Thus from, $\left(
3.1\right) $ we can write%
\begin{equation}
x_{u}\times x_{v}\mid _{p}=\frac{2\left\langle \left( p-c\right)
,x_{u}\times x_{v}\mid _{p}\right\rangle }{\left\Vert p-c\right\Vert ^{2}}%
\left( p-c\right) .  \tag{3.2}
\end{equation}%
Taking the inner product of $\left( 3.2\right) $ with $p-c,$ we obtain 
\begin{equation*}
\left\langle x_{u}\times x_{v}\mid _{p},p-c\right\rangle =0.
\end{equation*}%
If this last equation takes into account equality in $\left( 3.2\right) $,
we have%
\begin{equation*}
x_{u}\times x_{v}\mid _{p}=0.
\end{equation*}%
So, the assumption is incorrect, that is, $\Phi \left( p\right) =q$ is the
regular point of $y.$

Now, we assume that $p=x\left( u_{0},v_{0}\right) $ is a singular point of $%
x,$ then 
\begin{equation*}
x_{u}\times x_{v}\mid _{p}=0.
\end{equation*}%
Thus from the $\left( 3.1\right) ,$ we obtain that%
\begin{equation*}
y_{u}\times y_{v}\mid _{q}=0.
\end{equation*}%
This completes the proof.\newline

Let\textbf{\ }$I,$ $II$ be the first and second fundamental forms of $x$,
and let $\widetilde{I}$, $\widetilde{II}$ be these of $y.$ From $\cite{8},$
we have \textbf{\ }%
\begin{equation}
\widetilde{I}\circ \Phi _{\ast }=\lambda ^{2}I,  \tag{3.3}
\end{equation}%
\begin{equation}
\widetilde{II}\circ \Phi _{\ast }=-\lambda II-2\delta I,  \tag{3.4}
\end{equation}%
where $\lambda \left( u,v\right) =\tfrac{r^{2}}{\left\Vert x\left(
u,v\right) -c\right\Vert ^{2}}$ and $\delta \left( u,v\right) =\frac{%
2r^{2}\left\langle U_{x}(u,v),\left( x\left( u,v\right) -c\right)
\right\rangle }{\left\Vert x\left( u,v\right) -c\right\Vert ^{4}}.$ \newline

Thus, we have following results without proofs.\newline

\textbf{Corollary 3.2. }$\left( \cite{1}\right) $ Let $\Gamma _{jk}^{i}$ and 
$\widetilde{\Gamma }_{jk}^{i}$ $\left( i,j,k=1,2\right) $ be Christoffel
symbols of $x$ and $y,$ respectively. Then%
\begin{equation*}
\begin{array}{cc}
\widetilde{\Gamma }_{11}^{1}=\Gamma _{11}^{1}+\dfrac{\left[ EG-2F^{2}\right] 
\dfrac{\partial \lambda ^{2}}{\partial u}+FE\dfrac{\partial \lambda ^{2}}{%
\partial v}}{2\lambda ^{2}\left[ EG-F^{2}\right] },\text{ \ } & \widetilde{%
\Gamma }_{11}^{2}=\Gamma _{11}^{2}+\dfrac{EF\dfrac{\partial \lambda ^{2}}{%
\partial u}-E^{2}\dfrac{\partial \lambda ^{2}}{\partial v}}{2\lambda ^{2}%
\left[ EG-F^{2}\right] },\text{ \ \ \ \ \ \ \ \ \ \ \ \ \ \ }%
\end{array}%
\end{equation*}%
\begin{equation*}
\begin{array}{cc}
\widetilde{\Gamma }_{12}^{1}=\Gamma _{12}^{1}+\dfrac{GE\dfrac{\partial
\lambda ^{2}}{\partial v}-GF\dfrac{\partial \lambda ^{2}}{\partial u}}{%
2\lambda ^{2}\left[ EG-F^{2}\right] },\text{ \ \ \ \ \ \ \ \ \ \ \ \ } & 
\widetilde{\Gamma }_{12}^{2}=\Gamma _{12}^{2}+\dfrac{EG\dfrac{\partial
\lambda ^{2}}{\partial u}-EF\dfrac{\partial \lambda ^{2}}{\partial v}}{%
2\lambda ^{2}\left[ EG-F^{2}\right] },\text{ \ \ \ \ \ \ \ \ \ \ \ \ \ }%
\end{array}%
\end{equation*}%
\begin{equation*}
\begin{array}{cc}
\widetilde{\Gamma }_{22}^{1}=\Gamma _{22}^{1}+\dfrac{GF\dfrac{\partial
\lambda ^{2}}{\partial v}-G^{2}\dfrac{\partial \lambda ^{2}}{\partial u}}{%
2\lambda ^{2}\left[ EG-F^{2}\right] },\text{ \ \ \ \ \ \ \ \ \ \ \ \ \ } & 
\widetilde{\Gamma }_{22}^{2}=\Gamma _{22}^{2}+\dfrac{\left[ EG-2F^{2}\right] 
\dfrac{\partial \lambda ^{2}}{\partial v}+FG\dfrac{\partial \lambda ^{2}}{%
\partial u}}{2\lambda ^{2}\left[ EG-F^{2}\right] },\text{ \ }%
\end{array}%
\end{equation*}%
where $\lambda \left( u,v\right) =\tfrac{r^{2}}{\left\Vert x\left(
u,v\right) -c\right\Vert ^{2}}$ and $E,F,G$ $\ $are coefficients of the
first fundamental form of $x.$\newline

\textbf{Corrolary 3.3.} $\left( \cite{1}\right) $ Let $\Gamma _{jk}^{i}$ and 
$\widetilde{\Gamma }_{jk}^{i}$ $\left( i,j,k=1,2\right) $ be Christoffel
symbols of $x$ and $y,$ respectively. Suppose that $x$ is a principal patch
in $E^{3}$, then%
\begin{equation*}
\widetilde{\Gamma }_{11}^{1}=\Gamma _{11}^{1}+\left( 2\lambda \right) ^{-2}%
\frac{\partial \lambda ^{2}}{\partial u},\text{ \ \ \ \ \ }\widetilde{\Gamma 
}_{11}^{2}=\Gamma _{11}^{2}-\left( 2\lambda \right) ^{-2}\frac{E}{G}\frac{%
\partial \lambda ^{2}}{\partial v},
\end{equation*}%
\begin{equation*}
\widetilde{\Gamma }_{12}^{1}=\Gamma _{12}^{1}+\left( 2\lambda \right) ^{-2}%
\frac{\partial \lambda ^{2}}{\partial v},\text{ \ \ \ \ \ }\widetilde{\Gamma 
}_{12}^{2}=\Gamma _{12}^{2}+\left( 2\lambda \right) ^{-2}\frac{\partial
\lambda ^{2}}{\partial u}\text{ \ \ \ \ }
\end{equation*}%
\begin{equation*}
\widetilde{\Gamma }_{22}^{1}=\Gamma _{22}^{1}-\left( 2\lambda \right) ^{-2}%
\frac{G}{E}\frac{\partial \lambda ^{2}}{\partial u},\text{ \ \ }\widetilde{%
\Gamma }_{22}^{2}=\Gamma _{22}^{2}+\left( 2\lambda \right) ^{-2}\frac{%
\partial \lambda ^{2}}{\partial v}\text{ \ \ \ \ }
\end{equation*}%
where $\lambda \left( u,v\right) =\tfrac{r^{2}}{\left\Vert x\left(
u,v\right) -c\right\Vert ^{2}}$ and $E,G$ $\ $are coefficients of the first
fundamental form of $x.$\newline

\textbf{Corollary 3.4 }$\left( \cite{1}\right) $ Let $S$ and $\tilde{S}$ \
be the shape operators of $x$ and $y$, respectively. Then%
\begin{equation}
\tilde{S}\circ \Phi _{\ast }=-\lambda ^{-1}S-\frac{2}{r^{2}}\eta I_{2}, 
\tag{3.8}
\end{equation}%
\newline
where $I_{2}$ is identity, $\lambda \left( u,v\right) =\tfrac{r^{2}}{%
\left\Vert x\left( u,v\right) -c\right\Vert ^{2}}$ and $\eta =\left\langle
U_{x},\left( x-c\right) \right\rangle $.\newline

\textbf{Corollary 3.5. }$\left( \cite{1}\right) $\textbf{\ (i) }The
inversion maps flat points into umbilic points,\newline
\textbf{(ii) }The inversion maps flat points into flat points if and only if 
\begin{equation*}
\eta =\left\langle U_{x},\left( x-c\right) \right\rangle =0.
\end{equation*}

\textbf{4.} \textbf{The} \textbf{Aymptotic and Conjugate} \textbf{Vectors in
Inverse Surfaces}\newline

\textbf{Lemma 4.1.} Denote by $III$ and $\widetilde{III}$ the third
fundamental forms of $x$ and $y,$ respectively. Then,%
\begin{equation}
\widetilde{III}\circ \Phi _{\ast }=III+2\xi II+\xi ^{2}I,  \tag{4.1}
\end{equation}%
where $\xi =\frac{2\left\langle U_{x},\left( x-c\right) \right\rangle }{%
\left\Vert x-c\right\Vert ^{2}}.$\newline

\textbf{Proof.\ \ }Assume that $p$ be any point in $x\left( U\right) $ and $%
\Phi \left( p\right) =q\in y\left( U\right) $. Let the tangent map of $\Phi $
at $p$ be $\Phi _{\ast p}=T$ $_{p}\left( x\right) \longrightarrow T_{\Phi
\left( p\right) }\left( y\right) .$

Suppose $u_{p},v_{p}\in T$ $_{p}\left( x\right) $ and $\widetilde{u}_{q},%
\widetilde{v}_{q}\in T_{q}\left( y\right) $ such that 
\begin{equation*}
\Phi _{\ast p}\left( u_{p}\right) =\widetilde{u}_{q}\text{ \ \ \ and \ \ \ }%
\Phi _{\ast p}\left( v_{p}\right) =\widetilde{v}_{q}.
\end{equation*}%
Then%
\begin{equation*}
\widetilde{III}_{q}\left( \widetilde{u}_{q},\widetilde{v}_{q}\right)
=\left\langle \tilde{S}^{2}\left( \widetilde{u}_{q}\right) ,\widetilde{v}%
_{q}\right\rangle _{q},
\end{equation*}

or%
\begin{equation}
\widetilde{III}_{q}\left( \widetilde{u}_{q},\widetilde{v}_{q}\right)
=\left\langle \tilde{S}\left( \widetilde{u}_{q}\right) ,\tilde{S}\left( 
\widetilde{v}_{q}\right) \right\rangle _{q}.  \tag{4.2}
\end{equation}%
Considering $(3.3)$ and $(3.8)$ in $(4.2)$, we write%
\begin{equation*}
\widetilde{III}_{q}\left( \widetilde{u}_{q},\widetilde{v}_{q}\right)
=\lambda ^{2}\left\langle \left( -\lambda ^{-1}S-\frac{2}{r^{2}}\eta
I_{2}\right) \left( u_{p}\right) ,\left( -\lambda ^{-1}S-\frac{2}{r^{2}}\eta
I_{2}\right) \left( v_{p}\right) \right\rangle _{p}.
\end{equation*}%
From above equation, we obtain.%
\begin{equation*}
\widetilde{III}\circ \Phi _{\ast }=III+2\xi II+\xi ^{2}I.
\end{equation*}

\textbf{Lemma 4.2. }Let $\Phi _{\ast p}=T$ $_{p}\left( x\right)
\longrightarrow T_{\Phi \left( p\right) }\left( y\right) .$ Then, the
conjugate vectors in $x$ are invariant under $\Phi _{\ast }$ if and only if 
\begin{equation*}
\eta =\left\langle U_{x},\left( x-c\right) \right\rangle =0
\end{equation*}
or these are perpendicular.\newline

\textbf{Proof.\ }Let $p\in x\left( U\right) $ and $\Phi \left( p\right)
=q\in y\left( U\right) .$ We assume that 
\begin{equation*}
\Phi _{\ast }\left( u_{p}\right) =\widetilde{u}_{q}\text{ \ \ \ and \ \ \ }%
\Phi _{\ast }\left( v_{p}\right) =\widetilde{v}_{q},
\end{equation*}%
where $u_{p},v_{p}\in T$ $_{p}\left( x\right) $ and $\widetilde{u}_{q},%
\widetilde{v}_{q}\in T_{q}\left( y\right) .$ Let $u_{p}$ and $v_{p}\in T$ $%
_{p}\left( x\right) $ be conjugate vectors. Then, 
\begin{equation}
II_{p}\left( u_{p},v_{p}\right) =0.  \tag{4.3}
\end{equation}%
By $\left( 3.4\right) ,$ we write that 
\begin{equation}
\widetilde{II}_{q}\left( \widetilde{u}_{q},\widetilde{v}_{q}\right) =\dfrac{%
-2r^{2}\left\langle U_{x}\mid _{p},\left( p-c\right) \right\rangle }{%
\left\Vert p-c\right\Vert ^{4}}I_{p}\left( u_{p},v_{p}\right) ,  \tag{4.4}
\end{equation}%
where $\left\langle U_{x}\mid _{p},\left( p-c\right) \right\rangle =\eta
\left( p\right) .$ If $\widetilde{u}_{q}$ and $\widetilde{v}_{q}$ are
conjugate vectors then

\begin{equation*}
\widetilde{II}_{q}\left( \widetilde{u}_{q},\widetilde{v}_{q}\right) =0.
\end{equation*}%
Thus from $\left( 4.4\right) ,$ we obtain that%
\begin{equation*}
\eta \left( p\right) =\left\langle U_{x}\mid _{p},\left( p-c\right)
\right\rangle =0,
\end{equation*}%
or%
\begin{equation*}
I_{p}\left( u_{p},v_{p}\right) =\left\langle u_{p},v_{p}\right\rangle =0.
\end{equation*}%
\qquad

Conversely, let $\eta \left( p\right) =0$ or let $u_{p}\perp v_{p}$. Then,
from $\left( 4.1\right) ,$ we write that 
\begin{equation*}
\widetilde{II}_{q}\left( \widetilde{u}_{q},\widetilde{v}_{q}\right) =0.
\end{equation*}%
This implies that $\widetilde{u}_{q}$ and $\widetilde{v}_{q}$ are conjugate
vectors. This completes the proof.\newline

\textbf{Theorem 4.3. }Let $III$ and $\widetilde{III}$ be the third
fundamental forms of $x$ and $y,$ respectively. Assume that $II\neq -\frac{%
\xi }{2}I$ and the conjugate vectors in $x$ are not perpendicular each
other, then following conditions are equivalent:\newline

\textbf{(i) }the\textbf{\ }conjugate vectors in $x$ are invariant under $%
\Phi _{\ast }$

\textbf{(ii) }$\widetilde{III}\circ \Phi _{\ast }=III$,

\textbf{(iii) }$U_{y}=-U_{x}$, where $U_{x}$ and $U_{y}$ are unit normals of 
$x$ and $y,$ respectively.\newline

Proof. We show that $(i)\Rightarrow \left( ii\right) \Rightarrow \left(
iii\right) \Rightarrow \left( i\right) .$ If $\left( i\right) $ holds then 
\begin{equation*}
\eta =\left\langle U_{x},\left( x-c\right) \right\rangle =0,
\end{equation*}%
thus $\xi =0.$ By $\left( 4.1\right) ,$ we obtain 
\begin{equation*}
\widetilde{III}\circ \Phi _{\ast }=III.
\end{equation*}

If $\left( ii\right) $ holds then 
\begin{equation*}
\xi =\frac{2\left\langle U_{x},\left( x-c\right) \right\rangle }{\left\Vert
x-c\right\Vert ^{2}}=0,
\end{equation*}
thus $\eta =0.$ Also, from $\cite{8},$ we write

\begin{equation}
U_{y}=-U_{x}+\frac{2\eta }{\left\Vert x-c\right\Vert ^{2}}\left( x-c\right) .
\tag{4.5}
\end{equation}%
Hence, by $\left( 4.5\right) ,$ we obtain%
\begin{equation*}
U_{y}=-U_{x}.
\end{equation*}

Conversely, if $\left( iii\right) $ holds then by $\left( 4.5\right) $%
\begin{equation*}
\eta =0,
\end{equation*}%
thus it follows from Lemma 4.2 that $\left( i\right) $ holds. This completes
the proof.\newline

\textbf{Lemma 4.4. }Let\textbf{\ }$k$ and $\widetilde{k}$ be normal
curvatures of $x$ and $y$, respectively.\ Then%
\begin{equation}
\widetilde{k}\circ \Phi _{\ast }=-\frac{1}{\lambda }k-\frac{2}{r^{2}}\eta 
\tag{4.6}
\end{equation}%
where $\lambda \left( u,v\right) =\frac{r^{2}}{\left\Vert x\left( u,v\right)
-c\right\Vert ^{2}}$ and $\eta \left( u,v\right) =\left\langle
U_{x}(u,v),\left( x\left( u,v\right) -c\right) \right\rangle .$\newline

\textbf{Proof.\ }Let $p\in x\left( U\right) ,$ $\Phi \left( p\right) =q\in
y\left( U\right) $ and $\Phi _{\ast p}=T$ $_{p}\left( x\right)
\longrightarrow T_{\Phi \left( p\right) }\left( y\right) .$ Suppose that 
\begin{equation*}
\Phi _{\ast }\left( v_{p}\right) =\widetilde{v}_{q},
\end{equation*}%
where $v_{p}\in T_{p}\left( x\right) $ and $\widetilde{v}_{q}\in T_{q}\left(
x\right) .$ The normal curvature of $y$ in the direction $\widetilde{v}_{q}$
is%
\begin{equation}
\widetilde{k}\left( \widetilde{v}_{q}\right) =\frac{\widetilde{II}_{q}}{%
\widetilde{I}_{q}}\left( \widetilde{v}_{q},\widetilde{v}_{q}\right) . 
\tag{4.7}
\end{equation}

On the other hand, from Lemma 2.2, we write 
\begin{equation}
\left\Vert \widetilde{v}_{q}\right\Vert =\lambda \left( p\right) \left\Vert
v_{p}\right\Vert .  \tag{4.8}
\end{equation}%
Thus, by $(4.7)$ and $(4.8)$, we obtain 
\begin{equation*}
\widetilde{k}\circ \Phi _{\ast }\left( v_{p}\right) =-\frac{1}{\lambda
\left( p\right) }k\left( v_{p}\right) -\frac{2}{r^{2}}\eta \left( p\right) .
\end{equation*}

\textbf{Lemma 4.5. }Let $\Phi _{\ast p}=T_{p}\left( x\right) \longrightarrow
T_{\Phi \left( p\right) }\left( y\right) .$ Then, the asymptotic vectors in $%
x$ are invariant under $\Phi _{\ast }$ if and only if 
\begin{equation*}
\eta =\left\langle U_{x},\left( x-c\right) \right\rangle =0.\newline
\end{equation*}

\textbf{Proof.\ }Let $p\in x,$ $\Phi \left( p\right) =q\in y$ and 
\begin{equation*}
\Phi _{\ast }\left( v_{p}\right) =\widetilde{v}_{q},
\end{equation*}%
where $v_{p}\in T$ $_{p}\left( x\right) $ and $\widetilde{v}_{q}\in
T_{q}\left( y\right) .$\textbf{\ }Suppose that $v_{p}\in T$ $_{p}\left(
x\right) $ is a asymptotic vector, then 
\begin{equation*}
k\left( v_{p}\right) =0,
\end{equation*}%
where $k$ is the normal curvature of $x.$ By Lemma 4.4, we can write%
\begin{equation}
\widetilde{k}\left( \widetilde{v}_{q}\right) =-\frac{2}{r^{2}}\eta \left(
p\right) .  \tag{4.9}
\end{equation}%
\qquad 

If asymptotic vectors in $x$ are invariant under $\Phi _{\ast },$ then also $%
\widetilde{v}_{q}$ is asymptotic and%
\begin{equation*}
\eta \left( p\right) =0.
\end{equation*}%
\qquad \qquad\ 

Conversely, if $\eta \left( p\right) =0$ then, from $(4.9)$ we obtain%
\begin{equation}
\widetilde{k}\left( \widetilde{v}_{q}\right) =0.  \tag{4.10}
\end{equation}%
This is needed.\newline

\textbf{Theorem 4.6. }Let $\left[ III\right] $ and $\left[ \widetilde{III}%
\right] $ be the third fundamental forms of $x$ and $y,$ respectively. Then
following conditions are equivalent:\newline

\textbf{(i) }the\textbf{\ }asymptotic vectors in $x$ are invariant under $%
\Phi _{\ast }$

\textbf{(ii) }$\left[ \widetilde{III}\right] \circ \Phi _{\ast }=\left[ III%
\right] $,

\textbf{(iii) }$U_{y}=-U_{x}$, where $U_{x}$ and $U_{y}$ are unit normals of 
$x$ and $y,$ respectively.\newline

\textbf{Proof.} The proof is the same with that of Theorem 4.3.\newline

If the velocity vector of a curve is asymptotic vector, then the
curve\linebreak\ is asymptotic. Thus, we have a similar of Theorem 4.6 for
asymptotic curves.\newline

Now, we define the inverse curves of the parameter curves of a surface.%
\newline

\textbf{Definition 2.7. }Let $x$ be a patch in $\mathbb{E}^{n^{\ast }}$ and
let $y$ be the inverse patch of $x$ with respect to $\Phi .$ Then, the
curves 
\begin{equation*}
u\longrightarrow y\left( u,v_{0}\right) =c+\frac{r^{2}}{\left\Vert x\left(
u,v_{0}\right) -c\right\Vert ^{2}}\left( x\left( u,v_{0}\right) -c\right)
\end{equation*}%
and%
\begin{equation*}
v\longrightarrow y\left( u_{0},v\right) =c+\frac{r^{2}}{\left\Vert x\left(
u_{0},v\right) -c\right\Vert ^{2}}\left( x\left( u_{0},v\right) -c\right)
\end{equation*}%
are called inverse parameter curves of $u-$ parameter and $v-$ parameter
curves of $x,$ respectively, where $\left( u_{0},v_{0}\right) $ is a
constant point of $U.$\newline

If parameter curves of a patch are asymptotic curves, then it is called
\linebreak asymptotic patch. Hence, a similar of theorem 4.6 for asymptotic
patchs can be given.\newline

\textbf{Theorem 4.8. }Let\textbf{\ }$\tau $ and $\tilde{\tau}$ be torsions
of the asymptotic curves on $x$ and $y$, respectively$.$ If 
\begin{equation*}
\eta =\left\langle U_{x},\left( x-c\right) \right\rangle =0
\end{equation*}%
then

\begin{equation*}
\tilde{\tau}=\pm \frac{1}{\lambda }\tau ,
\end{equation*}%
where $\lambda \left( u,v\right) =\frac{r^{2}}{\left\Vert x\left( u,v\right)
-c\right\Vert ^{2}}.$

\textbf{Proof.\ }If $K$ and $\tilde{K}$ respectively are Gauss curvatures of 
$x$ and $y$, then from $\cite{8},$ we write 
\begin{equation*}
\tilde{K}=\dfrac{1}{\lambda ^{2}}K+\dfrac{4}{r^{2}}\lambda ^{-1}\eta H+%
\dfrac{4}{r^{4}}\eta ^{2},
\end{equation*}%
where $H$ is the mean curvature of $x.$ Because $\eta =0,$ we have%
\begin{equation*}
\tilde{K}=\dfrac{1}{\lambda ^{2}}K.
\end{equation*}%
Since $\tau $ and $\tilde{\tau}$ are the torsions of the asymptotic curves
in $x$ and $y$, respectively, we obtain 
\begin{equation*}
\tilde{\tau}=\pm \frac{1}{\lambda }\tau .
\end{equation*}

\textbf{5. Conclusions}\newline

For inverse surfaces, the minimality, the second and third fundamental
forms, the shape operators, the normal curvatures, the conjugate and the
asymptoic vectors and the torsions of the asymptotic curves are invariant\
properties under the inversion when $\eta =0$. We comment the function $\eta 
$ to be identically zero as follows: If 
\begin{equation*}
\eta \left( u,v\right) =\left\langle U_{x}(u,v),\left( x\left( u,v\right)
-c\right) \right\rangle =0,
\end{equation*}%
then the tangent planes at all the points of the surface $x(u,v)$ pass
through the center of the sphere of inversion. In such a case, the
properties mentioned above are invariant under the inversion.\newline

\end{document}